\documentclass[a4paper,11pt]{article}

\usepackage{amsmath,amsfonts,amssymb}

\newcommand{\prs}{\langle\;,\;\rangle}
\newcommand{\too}{\longrightarrow}
\newcommand{\om}{\omega}
\newcommand{\esp}{\quad\mbox{and}\quad}

\newcommand{\G}{{\mathfrak{g}}}

\newcommand{\ad}{{\mathrm{ad}}}
\newcommand{\tr}{{\mathrm{tr}}}

\newcommand{\B}{{\cal B}}

\newcommand{\al}{\alpha}
\newcommand{\be}{\beta}
\newcommand{\ga}{\gamma}

\newcommand{\la}{\lambda}

\font\bb=msbm10

\def\B{\hbox{\bb B}}
\def\R{\hbox{\bb R}}

\def\N{\hbox{\bb N}}

\newtheorem{theo}{Theorem}[section]
\newtheorem{pr}{Proposition}[section]
\newtheorem{Le}{Lemma}[section]

\newtheorem{rem}{Remark}

\title{Left-invariant  Lorentzian flat metrics on  Lie groups }          
\author{Malika Aitbenhaddou-Mohamed Boucetta-Hicham Lebzioui}                 

\begin{document}
 \maketitle

\begin{abstract} We call  the Lie algebra of a Lie group with a left invariant  
pseudo-Riemannian flat metric 
pseudo-Riemannian flat Lie algebra. We give a new proof of a classical result of Milnor on
Riemannian flat Lie algebras. We reduce the study of Lorentzian flat Lie algebras to those with
trivial center or those with  degenerate center. We show that the
double extension
process
can be used to construct all Lorentzian flat Lie algebras with  degenerate center
generalizing a result of Aubert-Medina on Lorentzian flat nilpotent Lie algebras. Finally, we give
the list of Lorentzian flat Lie algebras with  degenerate center
 up to
dimension 6.\\
\end{abstract}
{\it 2000 Mathematical Subject Classification: 53C50, 16T25; Secondary 53C20, 17B62.

}

\section{Introduction}\label{introduction}

A {\it pseudo-Riemannian flat Lie group} is a Lie group with a left invariant 
pseudo-Riemannian flat metric. The Lie algebra of such a Lie group is called
{\it pseudo-Riemannian  flat Lie algebra}. If the metric on the Lie group is complete the Lie
algebra is called  complete. It is a well-known result that a pseudo-Riemannian Lie algebra is
complete if and only if it is unimodular. A Riemannian (resp. Lorentzian) flat Lie group
is a pseudo-Riemannian flat Lie group for which the metric is definite positive (resp. of
signature $(-,+\ldots+)$). In \cite{Milnor}, Milnor showed that a Lie group is a Riemannian flat
Lie group  if and only
if its Lie algebra is a semi-direct product of an abelian algebra $\mathfrak{b}$ with an
abelian ideal $\mathfrak{u}$ and, for any $u\in\mathfrak{b}$, $\ad_u$ is skew-symmetric.
The characterization of Lorentzian flat Lie algebras (eventually complete) is an open problem. It is
a well-known result  that a Lorentzian flat Lie algebra must be solvable (see \cite{della}). On the
other hand,  in \cite{Aub-Med}, Aubert and Medina showed that Nilpotent Lorentzian flat
Lie algebras are obtained by a double extension process from Riemannian abelian Lie
algebras. In this paper, we reduce the problem of finding Lorentzian flat Lie algebras to
the determination of Lorentzian flat Lie algebras with  degenerate center and
those with trivial center. We show that the
double extension
process
can be used  to construct all Lorentzian flat Lie algebras with  degenerate
center from  Riemannian flat Lie
algebras and the Lie algebras obtained are unimodular and hence complete (see Theorem \ref{main}).
This result
generalizes Aubert-Medina's result.
We give
the list of Lorentzian flat Lie algebras with  degenerate center
 up to
dimension 6. The paper is organized as follows.  In Section
\ref{preliminaries}, we recall the
double extension process and state our main result (Theorem \ref{main}). In Section
\ref{proof}, we revisit Milnor's result and give a new formulation and a new proof of this
theorem  using the Lie
algebra of left invariant Killing vector fields of
a pseudo-Riemannian flat Lie group (See Theorem \ref{milnor}). This Lie algebra will play a crucial
role in the proof of our
main
result. Indeed, we will
establish a key Lemma (See Lemma \ref{le1}) involving this Lie algebra and, as a consequence, we
reduce the problem
of finding Lorentzian flat Lie algebras to
the determination of Lorentzian flat Lie algebras with  degenerate center and
those with trivial center, we recover the  keystone in the proof of Aubert-Medina's result (see
\cite{Aub-Med} Lemma 1.1) and prove Theorem \ref{main}. In Section \ref{example}, we give some
indications on how one can construct the tools used in the double extension process and we give
the list of Lorentzian flat Lie algebras with  degenerate center
 up to
dimension 6.\\

\section{Statement of the main result}\label{preliminaries}
A Lie group $G$ together with a left-invariant pseudo-Riemannian metric  is called a
\emph{pseudo-Riemannian Lie group}. 
The left-invariant pseudo-Riemannian metric 
defines an inner product $\prs$ on the Lie algebra $\G$ of $G$, and conversely, any inner
product on $\G$ gives rise to an unique  left-invariant metric on $\G$. The couple
$(\G,\prs)$ is called {\it pseudo-Riemannian Lie algebra}. We use the adjective Riemannian
(resp. Lorentzian)
instead of pseudo-Riemannian when the metric is definite positive (resp. of signature
$(-,+\ldots+)$). For any endomorphism $D:\G\too\G$ we denote
by $D^*:\G\too\G$ its adjoint with respect to $\prs$.\\
Let $(\G,\prs)$ be a pseudo-Riemannian Lie algebra of dimension $n$. The Levi-Civita
connection defines a product $(u,v)\mapsto uv$ on $\G$ called \emph{Levi-Civita
product} given by the Koszul formula
\begin{equation}\label{eq1}2\langle
uv,w\rangle=\langle[u,v],w\rangle+\langle[w,u],v\rangle+\langle[w,v],u\rangle.
\end{equation}
For any $u\in\G$, we denote by $\mathrm{L}_u:\G\too\G$ and $\mathrm{R}_u:\G\too\G$,
respectively, the left multiplication and the right multiplication by $u$ given by
$$\mathrm{L}_uv=uv\esp\mathrm{R}_uv=vu.$$We have
\begin{equation}\label{eq2}\ad_u=\mathrm{L}_u-\mathrm{R}_u,\end{equation}where
$\ad_u:\G\too\G$ is given by $\ad_uv=[u,v]$. The \emph{mean curvature vector} on $\G$ is the vector
given by
\begin{equation}\label{mean}\langle H,u\rangle=\mathrm{tr(\ad_u)},\;\forall u\in\G.\end{equation}The
Lie algebra $\G$ is unimodular if and only if $H=0$. 
The curvature of $\prs$ is given by
$$\mathrm{K}(u,v)=\mathrm{L}_{[u,v]}-[\mathrm{L}_u,\mathrm{L}_v].$$
$(\G,\prs)$ is called \emph{pseudo-Riemannian flat Lie algebra} if $\mathrm{K}$ vanishes
identically. This is equivalent to the fact that $\G$ endowed with the Levi-Civita product
is a left symmetric
algebra, i.e., for any $u,v,w\in\G$,
$$\mathrm{ass}(u,v,w)=\mathrm{ass}(v,u,w),$$where
$\mathrm{ass}(u,v,w)=(uv)w-u(vw).$ This relation is equivalent to
\begin{equation}\label{left}\mathrm{R}_{uv}-\mathrm{R}_v\circ\mathrm{R}_u=[\mathrm{L}_u,\mathrm{R}_v
],\end{equation}for any $u,v\in\G$.
\\ 

Let us recall now the double extension process and some related
results as elaborated in \cite{Aub-Med}. Let
$(B,[\;,\;]_0,\prs_0)$ be a pseudo-Riemannian flat Lie algebra, $\xi,D:B\too B$ two
endomorphisms of $B$, $b_0\in B$ and $\mu\in\R$ such that:
\begin{enumerate}
 \item \label{enu1}$\xi$ is a 1-cocycle of  $(B,[\;,\;]_0)$ with respect to the
representation $\mathrm{L}:B\too\mathrm{End}(B)$ defined by the left multiplication
associated to the Levi-Civita product, i.e., for any $a,b\in B$,
\begin{equation}\label{eq3}
 \xi([a,b])=\mathrm{L}_a\xi(b)-\mathrm{L}_b\xi(a),
\end{equation}
\item \label{enu2}$D$ is a derivation of $(B,[\;,\;]_0)$,
\item $D-\xi$ is skew-symmetric with respect to $\prs_0$, 
\begin{equation} \label{eq5}
[D,\xi]=\xi^2-\mu\xi-\mathrm{R}_{b_0},
\end{equation}and for any $a,b\in B$
\begin{equation} \label{eq6}
a\xi(b)-\xi(ab)=D(a)b+aD(b)-D(ab).
\end{equation}

\end{enumerate}
 We call  $(\xi,D,\mu,b_0)$ satisfying the three conditions above
\emph{admissible}.\\ Given $(\xi,D,\mu,b_0)$ admissible, 
we endow the vector space $\G=\R z\oplus B\oplus\R \bar{z}$  with the inner product
$\prs$ which extends $\prs_0$, for which $\mathrm{span}\{z,\bar{z}\}$ and $B$ are orthogonal,
 $\langle z,z\rangle=\langle \bar{z},\bar{z}\rangle=0$ and $\langle z,\bar{z}\rangle=1$. We define
also on
$\G$ the bracket 
\begin{equation}\label{bracket}[\bar{z},z]=\mu z,\; [\bar{z},a]=D(a)-\langle
b_0,a\rangle_0z\esp[a,b]=[a,b]_0+\langle(\xi-\xi^*)(a),b\rangle_0z,\end{equation}where $a,b\in B$
and
$\xi^*$ is the adjoint of $\xi$ with respect to $\prs_0$. Then $(\G,[\;,\;],\prs)$ is a
pseudo-Riemannian flat Lie algebra  called \emph{double extension}
of $(B,[\;,\;]_0,\prs_0)$ according to $(\xi,D,\mu,b_0)$. Moreover, any nilpotent Lorentzian flat
Lie algebra is a double extension of a
Riemannian abelian Lie algebra according to $(\xi,D,\mu,b_0)$ with $\xi=D$, $D^2=0$ and
$\mu=0$.\\
We can now state  our main result.
\begin{theo}\label{main}
 A Lorentzian Lie algebra  with  degenerate center is flat 
if and only if it is a  double extension of a
Riemannian  flat Lie algebra $(B,[\;,\;]_0,\prs_0)$ according to $(\xi,D,0,b_0)$ with
$(D,b_0)\not=(0,0)$. Moreover, a Lorentzian flat Lie algebra with degenerate center is unimodular
and hence complete.
\end{theo}

\section{The Lie algebra of left invariant Killing vector fields of
a pseudo-Riemannian flat Lie group}\label{proof}
\paragraph{The Riemannian case: Milnor theorem revisited}
In this paragraph, we give a new formulation and a new proof of Milnor's theorem  using the Lie
algebra of left invariant Killing vector fields of
a Riemannian flat Lie group. This Lie algebra will play a crucial role in the proof of our main
result.

 Let $(\G,[\;,\;],\prs)$ be a
pseudo-Riemannian Lie algebra. The Lie subalgebra
\begin{equation}\label{killing}\mathfrak{L}(\G)=\left\{u\in\G,\ad_u+\ad_u^*=0\right\}
=\left\{u\in\G,\mathrm{R}_u+\mathrm{R}_u^*=0\right\}\end{equation}is called  \emph{Killing}
subalgebra of $\G$. Indeed, if $\G$ is the Lie algebra of left invariant vector fields of
a pseudo-Riemannian Lie group then $\mathfrak{L}(\G)$ is the Lie algebra of left invariant
Killing vector fields. On the other hand,  one can see easily that the
orthogonal of
the derived ideal of $\G$ is given by
\begin{equation}
\label{eq7}\mathfrak{D}(\G)^\perp=\{u\in\G,\mathrm{R}_u=\mathrm{R}_u^*\}.\end{equation}
Finally, we put
$$N_\ell(\G)=\left\{u\in\G,\mathrm{L}_u=0\right\}\esp
N_r(\G)=\left\{u\in\G,\mathrm{R}_u=0\right\}.$$We have obviously
\begin{equation}\label{eq11}N_r(\G)=(\G\G)^\perp.\end{equation}
\begin{pr}\label{pr1}
 Let  $(\G,[\;,\;],\prs)$ be a
pseudo-Riemannian flat Lie algebra. Then:
\begin{enumerate}
 \item \label{pr11} For any $u\in \mathfrak{L}(\G)$, $\mathrm{R}_u^2=0$ and
$[\mathrm{R}_u,\mathrm{L}_u]=0.$
\item \label{pr12} For any $u\in\mathfrak{D}(\G)^\perp$,  $\mathrm{R}_u$ is nilpotent  and
$[\mathrm{R}_u,\mathrm{L}_u]=\mathrm{R}_u^2.$
\item The mean curvature vector satisfies $H\in \mathfrak{D}(\G)\cap\mathfrak{D}(\G)^\perp$. In
particular, if $\G$ is non unimodular then $\mathfrak{D}(\G)$ is degenerate.
\end{enumerate}

\end{pr}
{\bf Proof} By using \eqref{eq1} one can see easily that, for any
$u\in\mathfrak{L}(\G)\cup \mathfrak{D}(\G)^\perp$, $u.u=0$ and deduce from \eqref{left} that
$$[\mathrm{R}_u,\mathrm{L}_u]=\mathrm{R}_u^{2}.$$
If $u\in \mathfrak{L}(\G)$ then $\mathrm{R}_u$ is skew-symmetric and, since $\mathrm{L}_u$ is always
skew-symmetric,   $[\mathrm{R}_u,\mathrm{L}_u]$ is skew-symmetric. But $\mathrm{R}_u^{2}$
is symmetric which implies  1. \\
On the other hand, one can deduce by induction that for any $k\in\N^*$
$$[\mathrm{R}_u^k,\mathrm{L}_u]=k\mathrm{R}_u^{k+1}.$$and hence
$\tr(\mathrm{R}_u^{k})=0$
for any $k\geq2$ which implies that $\mathrm{R}_u$ is nilpotent.\\
Since, for any $u,v\in\G$, $\tr(\ad_{[u,v]})=0$, we deduce that $H\in\mathfrak{D}(\G)^\perp$.
Now, for any $u\in\mathfrak{D}(\G)^\perp$, $\mathrm{R}_u$ is nilpotent and hence
$$\tr(\ad_u)=\tr(\mathrm{R}_u)=\langle H,u\rangle=0,$$which implies
$H\in\mathfrak{D}(\G)$.\hfill$\square$

\begin{rem}
 If $\G$ is a Lorentzian flat Lie algebra, one can deduce from Proposition \ref{pr1} that for any
$u\in \mathfrak{D}(\G)^\perp$, $\mathrm{R}_u^3=0$. Moreover, if $\G$ is non unimodular then
$$\mathfrak{D}(\G)\cap\mathfrak{D}(\G)^\perp=\R H.$$
\end{rem}

We can give now a new formulation and a new proof of Milnor's theorem (see \cite{Milnor}). This new
formulation appeared first in \cite{bou}.
\begin{theo}\label{milnor} Let $G$ be a Riemannian Lie group. Then the curvature of $G$ vanishes
if and only if  $\mathfrak{L}(\G)$ is abelian, 
$\mathfrak{D}(\G)$ is abelian and $\mathfrak{L}(\G)^\perp=\mathfrak{D}(\G)$. Moreover, in this
case if $\dim\mathfrak{L}(\G)\geq1$ then the dimension of $\mathfrak{D}(\G)$ is even and the
Levi-Civita product is given
by\begin{equation}\label{eq12}\mathrm{L}_a=\left\{\begin{array}{ccc}\ad_a&\mbox{if}&a\in\mathfrak{L}
(\G),
\\ 0&\mbox{if}&a\in\mathfrak{D}(\G).
                      \end{array}\right.\end{equation}
\end{theo}
{\bf Proof.} Suppose that $G$ is a Riemannian flat Lie group. Since in definite positive context
skew-symmetric or
symmetric nilpotent endomorphism must vanish, we deduce from Proposition \ref{pr1} that
\begin{equation}\label{eq10}\mathfrak{L}(\G)=\mathfrak{D}(\G)^\perp=N_r(\G)=(\G\G)^\perp.
\end{equation}
These relations implies that $\mathfrak{L}(\G)$ is abelian, $\mathfrak{D}(\G)=\G\G$ and hence
$\mathfrak{D}(\G)$ is a bilateral ideal of the Levi-Civita product so $\mathfrak{D}(\G)$ endowed
with the restricted metric is a Riemannian flat  Lie algebra. It is known that a left
symmetric algebra cannot be equal to its derived ideal (see \cite{helms} pp.31). Hence  
$\G\not=\mathfrak{D}(\G)$,    and since $\mathfrak{D}(\G)$ is also flat
$\mathfrak{D}(\G)\not=\mathfrak{D}(\G)^2$ and so on. So $\G$ must
be solvable and
$$\mathfrak{D}(\G)=\mathfrak{L}(\mathfrak{D}(\G))\oplus \mathfrak{D}^2(\G).$$ Now the center of
$\mathfrak{D}(\G)$
is contained in $\mathfrak{L}(\mathfrak{D}(\G))$ and it intersects non trivially
$\mathfrak{D}^2(\G)$
($\mathfrak{D}(\G)$ is nilpotent) so $\mathfrak{D}(\G)$ must be abelian. This achieves the direct
part of the theorem. The equation \eqref{eq12} is easy to establish and the converse follows
immediately from this equation.\\
Suppose that $G$ is flat, and $\dim\mathfrak{L}(\G)\geq1$. Hence $\mathfrak{L}(\G)$ is abelian, 
$\mathfrak{D}(\G)$ is abelian and
$$\G=\mathfrak{L}(\G)\oplus\mathfrak{D}(\G).$$
Let $(s_1,...,s_p)$ be a basis of $\mathfrak{L}(\G)$.
The restriction of $\ad_{s_1}$ to $\mathfrak{D}(\G)$ is 
a skew-symmetric endomorphism, thus its kernel $K_1$ is of even codimension in $\mathfrak{D}(\G)$.
Now,
$\ad_{s_2}$ commutes with $\ad_{s_1}$ and  $K_1$ is invariant by $\ad_{s_2}$. By using the same
argument as above, we deduce
 that $K_1\cap\ker\ad_{s_2}$ is of even codimension in $K_1$. Finally $K_1$ is of even 
 codimension in $\mathfrak{D}(\G)$. Thus, by induction, we show that
$$K_p=\mathfrak{D}(\G)\cap\left(\cap_{i=1}^p\ker\ad_{s_i}\right)$$ is an even codimensional
subspace
of $\mathfrak{D}(\G)$. Now from its definition  $K_p$ is contained in the center of $\G$ which is
contained
in $\mathfrak{L}(\G)$ and then $K_p=\{0\}$ and the second part of the theorem follows.
\hfill $\square$

\paragraph{The Lorentzian case} It is known that  a left invariant affine structure on a Lie group
$G$ is complete if
and only if for any $u\in\G$, $\mathrm{R}_u$ is nilpotent (see \cite{segal} for instance). If $G$ is
a Riemannian flat Lie group then the underline left invariant affine structure is complete and one
can
deduce \eqref{eq10} immediately. We have avoided to use this argument in the proof of Theorem
\ref{milnor} and we have used arguments which are not specific to the Riemannian case.
Unfortunately,  in the
Lorentzian case the argument used
 to prove \eqref{eq10} cannot be used  since in the Lorentzian
context there is non trivial   skew-symmetric or symmetric nilpotent endomorphisms. However, in the
following
lemma, we show that a part of\eqref{eq10} is still valid in the Lorentzian case.

\begin{Le}\label{le1}
 Let $(\G,[\;,\;],\prs)$ be a
Lorentzian flat Lie algebra. Then
$$\mathfrak{L}(\G)=N_r(\G)=(\G\G)^\perp.$$
\end{Le}
{\bf Proof.}
Note first that we have always
 $N_r(\G)\subset\mathfrak{L}(\G)$.
 Let $u\in \mathfrak{L}(\G)$. According to Proposition \ref{pr1} \ref{pr11},
$\mathrm{R}_u^2=0$ and since $\mathrm{R}_u$ is skew-symmetric we get   that
$\mathrm{Im}\mathrm{R}_u$ is a totally
isotropic subspace and hence there exists an isotropic vector $e\in\G$ and a covector $\al\in\G^*$
such that $\mathrm{R}_u(u)=\al(u)e$ for any $u\in\G$. Choose a basis
$\{e,\bar{e},f_1,\ldots,f_{n-2}\}$ of $\G$ such that $\mathrm{span}\{e,\bar{e}\}$ and
$\mathrm{span}\{f_1,\ldots,f_{n-2}\}$ are orthogonal, $\{f_1,\ldots,f_{n-2}\}$ is orthonormal,
$\bar{e}$ is isotropic and $\langle e,\bar{e}\rangle=1$. We have, for any $i=1,\ldots,n-2$,
\begin{eqnarray*}
 \langle \mathrm{R}_u(e),\bar{e}\rangle&=&\al(e)=-\langle e,\mathrm{R}_u(\bar{e})\rangle=0,\\
 \langle \mathrm{R}_u(\bar{e}),\bar{e}\rangle&=&0=\al(\bar{e}),\\
 \langle \mathrm{R}_u(f_i),\bar{e}\rangle&=&\al(f_i)=-\langle f_i,\mathrm{R}_u(\bar{e})\rangle=0,
\end{eqnarray*} hence $\al=0$ and then $u\in N_r(\G)$ which achieves the proof of the lemma.
\hfill $\square$

 Let $(\G,[\;,\;],\prs)$ be a
Lorentzian flat Lie algebra. There are some interesting consequences of Lemma \ref{le1}:
\begin{enumerate}
 \item We have \begin{equation}\label{eq8} 
Z(\G)=N_\ell(\G)\cap N_r(\G)\subset\mathfrak{L}(\G)=(\G\G)^\perp\subset\mathfrak
{ D } (\G)^\perp.
 \end{equation}
 Thus $Z(\G)$ and $Z(\G)^\perp$ are ideals of $\G$ and are bilateral ideals of $\G$ endowed with the
Levi-Civita product, so if $Z(\G)$ is non degenerate we have
$$\G=Z(\G)\oplus Z(\G)^\perp$$ and both $Z(\G)$ and $Z(\G)^\perp$ are flat when endowed  with the
restricted metric.
\item If $Z(\G)$ is  degenerate, we have
\begin{equation}\label{eq9}
  Z(\G)\cap \mathfrak{D}(\G)\subset Z(\G)\cap Z(\G)^\perp\subset N_\ell(\G)\cap N_r(\G),
 \end{equation}and hence $Z(\G)\cap Z(\G)^\perp$ is one dimensional bilateral ideal and its
orthogonal is also a bilateral ideal, so we can use the double extension process.
\item If $\G$ is nilpotent then $Z(\G)\cap \mathfrak{D}(\G)\not=\{0\}$ and hence
 $Z(\G)\cap Z(\G)^\perp$ is one dimensional ideal contained in $N_\ell(\G)$ and Lemma 1.1 of
\cite{Aub-Med} follows. Note that this lemma played a crucial role in the proof of Aubert-Medina's
result and we give here a generalization and a new proof of this lemma.
\end{enumerate}
In conclusion, we have shown  that the problem
of finding Lorentzian flat Lie algebras reduces to
the determination of solvable Lorentzian flat Lie algebras with  degenerate center and
those with trivial center.
\paragraph{Proof of Theorem \ref{main}} Suppose that $(\G,[\;,\;],\prs)$ is a
Lorentzian flat Lie algebra and $Z(\G)$ is   degenerate. Then $\mathfrak{I}
=Z(\G)\cap Z(\G)^\perp=\R z$ where $z$ is an isotropic vector. Now from \eqref{eq8} we
deduce that $\mathrm{L}_z=\mathrm{R}_z=0$ and hence $\mathfrak{I}$ is a bilateral ideal
for the Levi-Civita product. Moreover, the orthogonal $\mathfrak{I}^\perp$ is also
a bilateral ideal.
So, according to Proposition 3.1 of \cite{Aub-Med},  
$\G=\R z\oplus B\oplus\R \bar{z}$ and it is  a double extension of $B$ according to
$(\xi,D,\mu,b_0)$. The Lie bracket is given by
$$[\bar{z},z]=\mu z,\; [\bar{z},a]=D(a)-\langle
b_0,a\rangle_0z\esp[a,b]=[a,b]_0+\langle(\xi-\xi^*)(a),b\rangle_0z.$$
 Since $z\in Z(\G)$ then $\mu=0$.
 The converse is obviously true.\\
 On the other hand, according to the brackets above and the fact that $B$ is unimodular, $\G$ is
unimodular if and only if $\tr(D)=0$. In Section  \ref{example}, we will show that if
$(\xi,D,0,b_0)$ is admissible and $B$ abelian then $D-\xi$ is skew-symmetric and $\xi$ is nilpotent
so $\tr(D)=0$. When $B$ is non abelian the relation $\tr(D)=0$ follows from Proposition \ref{prs}.
This achieves the proof.
 \hfill $\square$

\begin{rem} Let $\G$ be a Lorentzian flat non unimodular Lie algebra. According to Theorem
\ref{main}, $Z(\G)$ and non degenerate and hence $\G=Z(\G)\oplus Z(\G)^\perp$. Moreover, we have
seen that $Z(\G)^\perp$ is a bilateral ideal with respect to Levi-Civita product so the restriction
of the metric to $Z(\G)^\perp$ is Lorentzian and flat. Thus we reduce the study of Lorentzian flat
non unimodular Lie algebras to those with trivial center. Moreover, if $\G$ is a such  Lie algebra
then
$\mathfrak{D}(\G)\cap\mathfrak{D}(\G)^\perp=\R H.$
 
\end{rem}

\section{Lorentzian flat Lie algebras with  degenerate center
up to dimension six}\label{example}
According to Theorem \ref{main}, one can determine entirely all Lorentzian flat Lie algebras with
 degenerate center if one can find all admissible  $(\xi,D,0,b_0)$ on Riemannian
flat Lie algebras. In this section, we will give a general method to solve the equations
satisfied by admissible $(\xi,D,0,b_0)$ and we will use this method to give explicitly the solutions
on Riemannian flat Lie algebras of dimension 2, 3 or 4.  This will permit us to establish the list
of all Lorentzian flat Lie algebras with  degenerate center
up to dimension six.
\paragraph{The abelian case} Let $B$ be a Riemannian flat abelian Lie algebra of dimension $n$. One
can see easily that  a data $(\xi,D,0,b_0)$ is admissible if and only if $A=D-\xi$ is 
skew-symmetric  and
\begin{equation}\label{eq13} [A,\xi]=\xi^2.\end{equation}
Let $(A,\xi)$ be a solution of \eqref{eq13} with $A$ is skew-symmetric. One can deduce by induction 
that, for
any $k\in\N^*$,
\begin{equation}\label{eq14}[A,\xi^k]=k\xi^{k+1},\end{equation} and hence, for any
$k\geq2$, 
\begin{equation}\tr(\xi^k)=0.\end{equation}This implies that $\xi$ is nilpotent. Thus
there exists $q\leq\dim B$ such that
$$\{0\}\not=\ker\xi\subsetneqq\ker\xi^2\subsetneqq\ldots\subsetneqq\ker\xi^{q}=B.$$
 Then we  have the orthogonal splitting of $B$
\begin{equation}\label{splitting}
 B=\bigoplus_{k=0}^{q-1}F_k,
\end{equation}
where $F_0=\ker\xi$ and, for any $k=1,\ldots,q-1$,
$F_k=\ker\xi^{k+1}\cap\left(\ker\xi^k\right)^\perp$.
The key point is that \eqref{eq14} implies
that $A(\ker\xi^k)\subset\ker\xi^k$ for any $k\in\N$ and since $A$ is skew-symmetric,
$A(F_k)\subset F_k$. By using an orthonormal basis which respect to the splitting
\eqref{splitting}, the matrix of $A$ and $\xi$ are simple and one can solve \eqref{eq13}
easily. The following remarks can be used to simplify the computations when solving \eqref{eq13}.
\begin{rem}\label{rem2}Let $(A,\xi)$ be a solution of \eqref{eq13} with $A$ is
skew-symmetric.\begin{enumerate}
                        \item \label{rem21}If $q=\dim B$ then for any $k=0,\ldots,n-1$, $\dim F_k=1$
and hence
$A=0$ so $\xi^2=0$ and then $\dim B=2$. So if $\dim B\geq3$ then $q<\dim B$.
\item \label{rem22} One can deduce easily from \eqref{eq14} that  $\ker A\subset\ker\xi^2$.
\item \label{rem23} For any $1\leq k\leq q-1$, we have
$$\dim F_k\leq\dim\ker\xi^j,\; j=1,\ldots,k+1.$$
                       \end{enumerate}\end{rem}
We can now use what above to find all admissible $(\xi,D,0,b_0)$ when $\dim B\leq4$.
\begin{pr}\label{pr3} Let $B$ be  a Riemannian flat abelian Lie algebra. Then:
 \begin{enumerate}
  \item \label{pr31}If $\dim B=2$ then $(\xi,D,0,b_0)$ is admissible if and only if
$\xi=D=0$ or there exists an
orthonormal
basis $\{e_1,e_2\}$ of $B$  such that the matrices of $\xi$ and $D$ in this basis are
$$\left(M(\xi)=M(D)=\begin{pmatrix}0&a\\0&0\end{pmatrix},\; a\not=0\right)
\quad\mbox{or}\quad\left(\xi=0,\; M(D)=\begin{pmatrix}0&\la\\-\la&0\end{pmatrix},\;\la>0\right).$$
\item If $\dim B=3$ then $(\xi,D,0,b_0)$ is admissible if and only if $\xi=D=0$ or there exists an
orthonormal
basis $\{e_1,e_2,e_3\}$ of $B$  such that the matrices of $\xi$ and $D$ in this basis are{\small
$$\left(M(\xi)=M(D)=\begin{pmatrix}0&0&a\\0&0&b\\0&0&0\end{pmatrix},\;(a, b)\not=(0,0)\right)
\;\mbox{or}\;\left(\xi=0,\;
M(D)=\begin{pmatrix}0&\la&0\\-\la&0&0\\0&0&0\end{pmatrix},\;\la>0\right).$$}
\item If $\dim B=4$ then $(\xi,D,0,b_0)$ is admissible if and only if $\xi=D=0$ or there exists an
orthonormal
basis $\{e_1,e_2,e_3,e_4\}$ of $B$  such that the matrices of $\xi$ and $D$ in this basis have 
one
of the following forms:{\small
\begin{eqnarray*}
 (f_1)\quad M(\xi)&=&M(D)=\begin{pmatrix}0&0&a&b\\0&0&c&d\\0&0&0&0\\0&0&0&0\end{pmatrix},\quad
ad-bc\not=0,\\
 (f_2)\quad M(\xi)&=&M(D)=\begin{pmatrix}0&0&0&a\\0&0&0&b\\0&0&0&c\\0&0&0&0\end{pmatrix},\quad
(a,b,c)\not=(0,0,0),\\
(f_3)\quad M(D)&=&\begin{pmatrix}0&a&0&0\\-a&0&0&0\\0&0&0&b\\0&0&-b&0\end{pmatrix},\; \xi=0,\;
(a,b)\not=(0,0),\\
(f_4)\quad M(D)&=&\begin{pmatrix}0&a&b&c\\-a&0&-c&b\\0&0&0&a\\0&0&-a&0\end{pmatrix},\;
  M(\xi)=\begin{pmatrix}0&0&b&c\\0&0&-c&b\\0&0&0&0\\0&0&0&0\end{pmatrix},\; a\not=0,\;
(b,c)\not=(0,0),\\
(f_5)\quad M(D)&=&\begin{pmatrix}0&a&0&0\\-a&0&0&0\\0&0&0&b\\0&0&0&0\end{pmatrix},\;
  M(\xi)=\begin{pmatrix}0&0&0&0\\0&0&0&0\\0&0&0&b\\0&0&0&0\end{pmatrix}\; a\not=0,\;
b\not=0.
\end{eqnarray*}}
\end{enumerate}
\end{pr}
{\bf Proof.} Note first that $(\xi,D,0,b_0)$ is admissible if and only if $A=D-\xi$ is
skew-symmetric and $(A,\xi)$ is a solution of \eqref{eq13}. If it is the case,  we will use
repeatedly
the fact that $A$ leaves invariant the splitting \eqref{splitting} and Remark \ref{rem2}.   Let
$(\xi,D,0,b_0)$ be admissible such that $(\xi,D)\not=(0,0)$. We consider the
integer $q$ defined in \eqref{splitting}.\begin{enumerate}
              \item When $\dim B=2$, there is two possibilities of $q$. If  $q=2$ then, according to
\eqref{splitting}, $B=\ker\xi\oplus
F_1$ and hence $A=0$. Thus  $D=\xi$ and $\xi^2=0$. If $q=1$ then
$\xi=0$ and $D$ is skew-symmetric.
\item When $\dim B=3$ then, according to Remark \ref{rem2}, there is also two possibilities of $q$. 
 If $q=1$ then $\xi=0$ and  $D$ is skew-symmetric.\\
  If $q=2$ then, according to \eqref{splitting} and Remark \ref{rem2}, $B=\ker\xi\oplus F_1$ with
$\dim\ker\xi\geq \dim F_1$. Thus 
  $\dim\ker\xi=2$ and $\dim F_1=1$.
 So there exists an orthonormal basis of $B$ such that
  $$M(A)=\begin{pmatrix}0&a&0\\-a&0&0\\0&0&0\end{pmatrix}\esp M(\xi)
  =\begin{pmatrix}0&0&b\\0&0&c\\0&0&0\end{pmatrix},$$where $a\in\R$ and $(b,c)\not=(0,0)$.
  A direct computation shows that \eqref{eq13} is equivalent to $A=0$ and hence $D=\xi$.
  \item When $\dim B=4$ then,  according to Remark \ref{rem2},  $q\leq3$. 
\begin{itemize}
 \item If $q=1$  then $\xi=0$ and $D$ is skew-symmetric, this gives $(f_3)$.
 \item If $q=2$ then $\xi\not=0$, $B=\ker\xi\oplus F_1$ and $\dim\ker\xi\geq\dim F_1$. We
distinguish two cases.\\
 {\bf First case:} $\dim\ker\xi=2$ and $\dim F_1=2$. Then there exists an orthonormal basis such
that
 $$M(A)=\begin{pmatrix}0&a&0&0\\-a&0&0&0\\0&0&0&b\\0&0&-b&0\end{pmatrix}
 \esp M(\xi)=\begin{pmatrix}0&0&c&d\\0&0&e&f\\0&0&0&0\\0&0&0&0\end{pmatrix},$$where $a,b\in\R$ and
$cf-ed\not=0$.
 A direct computation shows that \eqref{eq13} is equivalent to
  $$\left\{\begin{array}{ccc}ae+bd&=&0,\\be+ad&=&0,\\af-bc&=&0,\\bf-ac&=&0.\end{array}\right.$$ 
 If $(a,b)=(0,0)$ we recover $(f_1)$. If $(a,b)\not=(0,0)$ then, since $\xi\not=0$, we get $a=b$ or
$a=-b$. If $a=b$  we recover $(f_4)$. If $a=-b$ we recover also $(f_4)$ when we permute  $e_3$
and $e_4$.\\
{\bf Second case:} $\dim\ker\xi=3$ and $\dim F_1=1$. Then there exists an orthonormal basis such
that
 $$M(A)=\begin{pmatrix}0&a&0&0\\-a&0&0&0\\0&0&0&0\\0&0&0&0\end{pmatrix}
 \esp M(\xi)=\begin{pmatrix}0&0&0&b\\0&0&0&c\\0&0&0&d\\0&0&0&0\end{pmatrix},$$where $a\in\R$ and
$(b,c,d)\not=(0,0,0)$.
 A direct computation shows that \eqref{eq13} is equivalent to $ab=0$ and $ac=0$. If $a=0$
we recover $(f_2)$
and if $b=c=0$ we recover $(f_5)$.
 \item If $q=3$ then $\xi^2\not=0$, $B=\ker\xi\oplus F_1\oplus F_2$ with $\dim
F_2\leq\dim(\ker\xi\oplus F_1)$ and
$\dim F_2\leq\dim\ker\xi$.
 Hence $\dim\ker\xi=2$ and $\dim\ker F_1=\dim\ker F_2=1$ then there exists an orthonormal basis such
that
 \begin{eqnarray*}
  M(A)&=&\begin{pmatrix}0&a&0&0\\-a&0&0&0\\0&0&0&0\\0&0&0&0\end{pmatrix},\;
 M(\xi)=\begin{pmatrix}0&0&b&c\\0&0&d&e\\0&0&0&f\\0&0&0&0\end{pmatrix}\esp\\
 M(\xi^2)&=&\begin{pmatrix}0&0&0&fb\\0&0&0&df\\0&0&0&0\\0&0&0&0\end{pmatrix},\end{eqnarray*}where
$f\not=0$.
 This case is impossible since $e_4\in\ker A$ and  $e_4\notin\ker\xi^2$ which is in contradiction
with Remark \ref{rem2}.\hfill $\square$
 
 \end{itemize}

             \end{enumerate}

\paragraph{The non abelian case} Let $B$ be a Riemannian flat non abelian Lie algebra of dimension
$n$. According to Theorem \ref{milnor},  $\mathfrak{L}(B)$ and
$\mathfrak{D}(B)$ are abelian and
$$B=\mathfrak{L}(B)\oplus\mathfrak{D}(B).$$Moreover,
$$\mathrm{L}_a=\left\{\begin{array}{ccc}\ad_a&\mbox{if}&a\in\mathfrak{L}(B),\\
                       0&\mbox{if}&a\in\mathfrak{D}(B).
                      \end{array}\right.$$
                      Since $\mathfrak{L}(B)$ is abelian and acts 
 on $\mathfrak{D}(B)$ by skew-symmetric endomorphisms, there exists a family of non
vanishing
vectors   $u_1,\ldots,u_r\in \mathfrak{L}(B)$  and an orthonormal basis $(f_1,\ldots,f_{2r})$
of
$\mathfrak{D}(B)$ such that, for any $j=1,\ldots,r$ and any $s\in \mathfrak{L}(B)$,
\begin{equation}\label{central}[s,f_{2j-1}]=\langle s,u_j\rangle
f_{2j}\quad\mbox{and}\quad
[s,f_{2j}]=-\langle s,u_j\rangle f_{2j-1}.\end{equation}
If $F$ is an endomorphism of $B$, we put, for any $u\in B$,  $F(u)=F_1(u)+F_2(u)$ where $F_1(u)\in
\mathfrak{L}(B)$ and $F_2(u)\in\mathfrak{D}(B)$, and we denote by
$\overline{F}_1\in\mathrm{End}(\mathfrak{L}(B))$ and
$\overline{F}_2\in\mathrm{End}(\mathfrak{D}(B))$, respectively, the restriction of $F_1$ to
$\mathfrak{L}(B)$ and the restriction of
$F_2$ to $\mathfrak{D}(B)$.
\begin{pr}\label{prs}With the notations and hypothesis above, $(\xi,D,0,b_0)$ is admissible if and
only if
$\overline{D}_1-\overline{\xi}_1$ and $\overline{D}_2-\overline{\xi}_2$ are skew-symmetric and, for
any $a,b\in \mathfrak{L}(B)$ and any $c\in\mathfrak{D}(B)$,
\begin{eqnarray}
 {D_1}_{|\mathfrak{D}(B)}&=&{\xi_1}_{|\mathfrak{D}(B)}=0,\quad(\xi_2-D_2)_{|\mathfrak{L}(B)}
=0,\label{eq1s}\\
     0&=& [D_2(a),b]+[a,D_2(b)],\label{eq2s}\\
          D_2([a,c])&=&[D_1(a),c]+[a,D_2(c)],\label{eq3s}\\
 \xi_2([a,c])&=&[a,\xi_2(c)],\label{eq4s}\\
\;[\overline{D}_1,\overline{\xi}_1]&=&\overline{\xi_1}^2,\label{eq5s}\\
\;[\overline{D}_2,\overline{\xi}_2]&=&\overline{\xi_2}^2,\label{eq6s}\\
\;[D_2,\xi_2](a)&=&\xi_2^2(a)+\xi_2\circ D_1(a)+[b_0,a].\label{eq7s}
    \end{eqnarray} Moreover, if $(\xi,D,0,b_0)$ is admissible then $\tr(D)=0$.
\end{pr}

{\bf Proof.} Recall that $(\xi,D,0,b_0)$ is admissible if and only if $D$ is a derivation of $B$,
$D-\xi$ is skew-symmetric and $(\xi,D,b_0)$ satisfy \eqref{eq3}-\eqref{eq6}.\\
Now $D$ is a derivation and $(\xi,D)$ satisfy \eqref{eq3} and \eqref{eq6} if and
only if, for any $a,b\in\mathfrak{L}(B)$ and any $c,d\in\mathfrak{D}(B)$,
    \begin{eqnarray*}
     0&=& [D_2(a),b]+[a,D_2(b)],\\
     0&=& [D_1(c),d]+[c,D_1(d)],\\
     D_2([a,c])&=&[D_1(a),c]+[a,D_2(c)],\quad
 D_1([a,c])=0,\\
\;[a,\xi_2(b)]&=&[b,\xi_2(a)],\\
 \xi_2([a,c])&=&[a,\xi_2(c)],\;\xi_1([a,c])=0,\\
\;[a,\xi_2(b)]&=&[a,D_2(b)],\quad\quad\quad(*)\\
\;[a,\xi(c)]-\xi([a,c])&=&[D(a),c]+[a,D(c)]-D([a,c]),\\
 \;[D_1(c),d]&=&0.\end{eqnarray*}
We get obviously  that ${D_1}_{|\mathfrak{D}(B)}={\xi_1}_{|\mathfrak{D}(B)}=0$. Moreover, from
$(*)$ we deduce that, for any $b\in\mathfrak{L}(B)$, $\xi_2(b)-D_2(b)$ is a central element and
since
the center of $B$ is contained in $\mathfrak{L}(B)$, we deduce that $(\xi_2-D_2)_{|\mathfrak{L}(B)}
=0$. On the other hand, one can see easily that if
${D_1}_{|\mathfrak{D}(B)}={\xi_1}_{|\mathfrak{D}(B)}=0$ and $(\xi_2-D_2)_{|\mathfrak{L}(B)}
=0$ then  $(\xi,D,0,b_0)$ satisfies \eqref{eq5} if and only if
$$[D_1,\xi_1]=\xi_1^2\esp [D_2,\xi_2]=\xi_2^2+\xi_2D_1-\mathrm{R}_{b_0}.$$
By evaluating the second equation respectively on $\mathfrak{L}(B)$ and $\mathfrak{D}(B)$ one can
conclude. Moreover, if $(\xi,D,0,b_0)$ is admissible then one can deduce easily from what above
that $\tr(D)=\tr(\overline{\xi}_1)+\tr(\overline{\xi}_2)=0$ ($\overline{\xi}_1$ and
$\overline{\xi}_2$ are nilpotent by virtue of \eqref{eq5s} and \eqref{eq6s}).
\hfill $\square$

Let us use this proposition two find $(D,\xi,0,b_0)$ admissible when $\dim B=3$ or $4$.
\begin{pr} Let $B$ be a Riemannian flat non abelian Lie algebra of dimension 3. Then there
exists an orthonormal basis $\{e_1,e_2,e_3\}$ of $B$ such that
$\mathfrak{L}(B)=\mathrm{span}\{e_1\}$, $\mathfrak{D}(B)=\mathrm{span}\{e_2,e_3\}$ and
$$[e_1,e_2]=\la e_3\esp[e_1,e_3]=-\la e_2,$$where $\la>0$. Moreover, $(D,\xi,0,b_0)$ is admissible
if and only if
$$M(\xi)=\begin{pmatrix}0&0&0\\a&0&0\\ b&0&0\end{pmatrix},\,
M(D)=\begin{pmatrix}0&0&0\\a&0&c\\ b&-c&0\end{pmatrix}
\esp b_0=b_1e_1+\frac{ca}\la e_2+\frac{cb}\la e_3.$$\end{pr}

{\bf Proof.} The existence of the orthonormal basis $\B=\{e_1,e_2,e_3\}$ in
which the Lie bracket is given by the relations above is a consequence of Theorem
\ref{milnor}.\\
Suppose now that  $(D,\xi,0,b_0)$ is admissible. Since $\mathfrak{L}(B)$ is a line, the
equation \eqref{eq2s} holds and we  deduce from \eqref{eq5s} and the fact
$\overline{\xi}_1-\overline{D}_1$ is skew-symmetric  that
$\overline{\xi}_1=\overline{D}_1=0$ and hence $\xi_1=D_1=0$. So $D_2$ and $\xi_2$ satisfy
the same equation, namely \eqref{eq4s}. This equation is equivalent to
\begin{eqnarray*}
 \xi_2([e_1,e_2])&=&[e_1,\xi_2(e_2)],\\
 \la\xi_2(e_3)&=&\la\langle \xi_2(e_2),e_2\rangle e_3-\la\langle \xi_2(e_2),e_3\rangle
e_2,\\
\xi_2([e_1,e_3])&=&[e_1,\xi_2(e_3)],\\
 -\la\xi_2(e_2)&=&\la\langle \xi_2(e_3),e_2\rangle e_3-\la\langle \xi_2(e_3),e_3\rangle
e_2.
\end{eqnarray*}
These equations are equivalent to
\begin{eqnarray*}
 \langle \xi_2(e_3),e_2\rangle&=&-\langle \xi_2(e_2),e_3\rangle\esp 
 \langle \xi_2(e_3),e_3\rangle=\langle \xi_2(e_2),e_2\rangle.
\end{eqnarray*}
Thus, since $D-\xi$ is skew-symmetric, the matrices of $D$ and $\xi$ in the basis $\B$
are of the following form
$$M(\xi)=\begin{pmatrix}0&0&0\\a&e&d\\b&-d&e\end{pmatrix}\esp
M(D)=\begin{pmatrix}0&0&0\\a&e&c\\b&-c&e\end{pmatrix}.$$
On the other hand, according to Proposition \ref{pr3} \ref{pr31}, the equation
\eqref{eq6s} and the fact that $\overline{\xi}_2-\overline{D}_2$ is skew-symmetric is
equivalent to 
$$\left(M(\overline{\xi}_2)=0\esp
M(\overline{D}_2)^*=-M(\overline{D}_2)\right)\;\mbox{or}\;
\left(M(\overline{\xi}_2)=
M(\overline{D}_2)\esp M(\overline{\xi}_2)^2=0\right).\eqno(*)$$
Or $M(\overline{\xi}_2)=\begin{pmatrix}e&d\\-d&e\end{pmatrix}$ and
$M(\overline{D}_2)=\begin{pmatrix}e&c\\-c&e\end{pmatrix}$. We deduce that $(*)$ is
equivalent to $e=d=0$. Now, the equation \eqref{eq7s} is equivalent to $b_0=
b_1e_1+\frac{ca}\la e_2+\frac{cb}\la e_3$. This achieves the proof.
\hfill $\square$

\begin{pr}\label{pr7}Let $B$ be a Riemannian flat non abelian Lie algebra of dimension 4.  Then
there
exists an orthonormal basis $\{e_1,e_2,f_1,f_2\}$ of $B$ such that
$\mathfrak{L}(B)=\mathrm{span}\{e_1,e_2\}$, $\mathfrak{D}(B)=\mathrm{span}\{f_1,f_2\}$ and
$$[e_i,f_1]=\la_i f_2\esp[e_i,f_2]=-\la_i f_1,$$where $\la_1,\la_2\geq0$ and
$(\la_1,\la_2)\not=(0,0).$ Moreover, $(\xi,D,0,b_0)$ is admissible if and only if
\begin{eqnarray*}
 M(D)&=&\begin{pmatrix}x\la_1\la_2&x\la_2^2&0&0\\-x\la_1^2&-x\la_1\la_2&0&0\\c\la_1&c\la_2&0&f\\
      d\la_1&d\la_2&-f&0\end{pmatrix},\\
 M(\xi)&=&\begin{pmatrix}
x\la_1\la_2&x\la_2^2&0&0\\-x\la_1^2&-x\la_1\la_2&0&0\\c\la_1&c\la_2&0&0\\
      d\la_1&d\la_2&0&0\end{pmatrix},\\
      b_0&=&b_1e_1+b_2e_2+fc f_1+fd f_2.
   \end{eqnarray*}

\end{pr}

 {\bf Proof.} The existence of the orthonormal basis $\B=\{e_1,e_2,f_1,f_2\}$ in
which the Lie bracket is given by the relations above is a consequence of Theorem
\ref{milnor}.\\
Suppose now that  $(D,\xi,0,b_0)$ is admissible.
The equation \eqref{eq2s} is equivalent to
\begin{eqnarray*}
    0&=&[D_2(e_1),e_2]+[e_1,D_2(e_2)],\\
    &=&-\la_2\langle D_2(e_1),f_1\rangle f_2+\la_2\langle D_2(e_1),f_2\rangle f_1
    +\la_1\langle D_2(e_2),f_1\rangle f_2-\la_1\langle D_2(e_2),f_2\rangle f_1\\
    &=&\langle \la_2D_2(e_1)-\la_1D_2(e_2),f_2\rangle f_1-\langle
\la_2D_2(e_1)-\la_1D_2(e_2),f_1\rangle f_2.
\end{eqnarray*}This is equivalent to
$\la_2D_2(e_1)=\la_1D_2(e_2)$ and since $(\la_1,\la_2)\not=(0,0)$ this is equivalent to the
existence of $(c,d)\in\R^2$ such that
$$D_2(e_1)=\la_1(cf_1+df_2)\esp D_2(e_2)=\la_2(cf_1+df_2).$$
The equation \eqref{eq3s} is equivalent to
\begin{eqnarray*}
\la_iD_2(f_2)&=&[D_1(e_i),f_1]+[e_i,D_2(f_1)]\\
&=&\la_1\langle D_1(e_i),e_1\rangle f_2+\la_2\langle D_1(e_i),e_2\rangle f_2+\la_i\langle
D_2(f_1),f_1\rangle f_2-\la_i\langle
D_2(f_1),f_2\rangle f_1,\\
-\la_iD_2(f_1)&=&[D_1(e_i),f_2]+[e_i,D_2(f_2)]\\
&=&-\la_1\langle D_1(e_i),e_1\rangle f_1-\la_2\langle D_1(e_i),e_2\rangle f_1+\la_i\langle
D_2(f_2),f_1\rangle f_2-\la_i\langle
D_2(f_2),f_2\rangle f_1.
\end{eqnarray*}
This is equivalent to
\begin{eqnarray*}
 \langle D_2(f_1),f_2\rangle&=&-\langle D_2(f_2),f_1\rangle,\\
\la_i\langle D_2(f_2),f_2\rangle&=&\la_1\langle D_1(e_i),e_1\rangle+\la_2\langle
D_1(e_i),e_2\rangle+\la_i\langle
D_2(f_1),f_1\rangle,\\
\la_i\langle D_2(f_1),f_1\rangle&=&\la_1\langle D_1(e_i),e_1\rangle+\la_2\langle
D_1(e_i),e_2\rangle+\la_i\langle
D_2(f_2),f_2\rangle.
\end{eqnarray*}
Which is equivalent to the existence of $(a,b)\in\R^2$ such that
\begin{eqnarray*}
 \langle D_2(f_1),f_2\rangle&=&-\langle D_2(f_2),f_1\rangle,\;
\langle D_2(f_2),f_2\rangle=\langle
D_2(f_1),f_1\rangle,\\
D_1(e_1)&=&a(\la_2e_1-\la_1e_2)\esp D_1(e_2)=b(\la_2e_1-\la_1e_2).
\end{eqnarray*}
The equation \eqref{eq4s} is equivalent to
\begin{eqnarray*}
\la_i\xi_2(f_2)&=&[e_i,\xi_2(f_1)]
=\la_i\langle  \xi_2(f_1),f_1\rangle f_2-\la_i\langle  \xi_2(f_1),f_2\rangle f_1,\\
-\la_i\xi_2(f_1)&=&[e_i,\xi_2(f_2)]
=\la_i\langle  \xi_2(f_2),f_1\rangle f_2-\la_i\langle  \xi_2(f_2),f_2\rangle f_1.
\end{eqnarray*}
This is equivalent to
\begin{eqnarray*}
 \langle \xi_2(f_2),f_1\rangle&=&-\langle  \xi_2(f_1),f_2\rangle\esp
 \langle \xi_2(f_2),f_2\rangle=\langle  \xi_2(f_1),f_1\rangle.
\end{eqnarray*}
In conclusion, and since $D-\xi$ is skew-symmetric, we get that the matrices of $D$ and $\xi$ in
the basis $\B$ have the following form
\begin{eqnarray*}
M( D)&=&\begin{pmatrix}a\la_2&b\la_2&0&0\\-a\la_1&-b\la_1&0&0\\c\la_1&c\la_2&\al&-\be\\
      d\la_1&d\la_2&\be&\al\end{pmatrix}\esp
 M(\xi)=\begin{pmatrix}
a\la_2&b\la_2+e&0&0\\-a\la_1-e&-b\la_1&0&0\\c\la_1&c\la_2&\al&-f-\be\\
      d\la_1&d\la_2&\be+f&\al\end{pmatrix}.
   \end{eqnarray*}
   On the other hand, according to Proposition \ref{pr3} \ref{pr31}, the equation
\eqref{eq6s} and the fact that $\overline{\xi}_2-\overline{D}_2$ is skew-symmetric is
equivalent to 
$$\left(M(\overline{\xi}_2)=0\esp
M(\overline{D}_2)^*=-M(\overline{D}_2)\right)\;\mbox{or}\;
\left(M(\overline{\xi}_2)=
M(\overline{D}_2)\esp M(\overline{\xi}_2)^2=0\right).\eqno(*)$$
Or $M(\overline{\xi}_2)=\begin{pmatrix}\al&-f-\be\\f+\be&\al\end{pmatrix}$ and
$M(\overline{D}_2)=\begin{pmatrix}\al&-\be\\\be&\al\end{pmatrix}$. We deduce that $(*)$ is
equivalent to $\al=f+\be=0$. In a similar way, the equation
\eqref{eq5s} and the fact that $\overline{\xi}_1-\overline{D}_1$ is skew-symmetric is
equivalent to 
$$\left(M(\overline{\xi}_1)=0\esp
M(\overline{D}_1)^*=-M(\overline{D}_1)\right)\;\mbox{or}\;
\left(M(\overline{\xi}_1)=
M(\overline{D}_1)\esp M(\overline{\xi}_1)^2=0\right).\eqno(**)$$
Or $M(\overline{\xi}_1)=\begin{pmatrix}a\la_2&b\la_2+e\\-a\la_1-e&-b\la_1\end{pmatrix}$ and
$M(\overline{D}_1)=\begin{pmatrix}a\la_2&b\la_2\\-a\la_1&-b\la_1\end{pmatrix}$. A careful checking
shows that $(**)$ is equivalent to $e=a\la_2-b\la_1=0$, i.e., $e=0$, $a=x\la_1$ and $b=x\la_2$. In
conclusion
\begin{eqnarray*}
M(
D)&=&\begin{pmatrix}x\la_1\la_2&x\la_2^2&0&0\\-x\la_1^2&-x\la_1\la_2&0&0\\c\la_1&c\la_2&0&f\\
      d\la_1&d\la_2&-f&0\end{pmatrix}\esp
 M(\xi)=\begin{pmatrix}
x\la_1\la_2&x\la_2^2&0&0\\-x\la_1^2&-x\la_1\la_2&0&0\\c\la_1&c\la_2&0&0\\
      d\la_1&d\la_2&0&0\end{pmatrix}.
   \end{eqnarray*}
Finally, the equation \eqref{eq7s} is equivalent to $b_0=b_1e_1+b_2e_2+fc f_1+fd f_2$.
\hfill $\square$

By using Theorem \ref{main} and Propositions \ref{pr3}-\ref{pr7}, let us give  the list of
Lorentzian flat Lie algebras with  degenerate center
up to dimension six. We proceed as follows:
\begin{enumerate}
 \item We pick an admissible solution $(\xi,D,0,b_0)$ found in Propositions \ref{pr3}-\ref{pr7},
 \item by using \eqref{bracket}, we compute the Lie brackets and we make an appropriate change of
basis to get a simple form of these brackets,
\item finally,  we give the matrix of the metric in the new basis
which we continue to denote by $\{z,\bar{z},e_1,\ldots,e_n\}$
($n=1,\ldots,4$).
\end{enumerate}

\begin{tabular}{|l|l|l|}
 \hline
 The non vanishing Lie brackets&The matrix of the metric& Nature\\
 \hline
$\begin{array}{l}
 [\bar{z},e_1]=z.
\end{array}$&{\small$
 \begin{pmatrix}
  0&\al&0\\\al &0&0\\0&0&1\end{pmatrix}$}&Heisenberg\\
 $\al>0$. &&\\
\hline\end{tabular}

\begin{tabular}{|l|l|l|}
 \hline
$\begin{array}{l}
 [\bar{z},e_1]=z.
\end{array}$&{\small$
 \begin{pmatrix}
  0&\al&0&0\\\al&0&0&0\\0&0&1&a\\0&0&a&1+a^2\end{pmatrix}$}&2-nilpotent\\
  $a\in\R$, $\al>0$.&&\\
\hline
$\begin{array}{l}
 [\bar{z},e_1]=az,\;[\bar{z},e_2]=e_1,\\\;[e_1,e_2]=-z.
\end{array}$&{\small$
 \begin{pmatrix}
  0&\al&0&0\\\al&0&\be&0\\0&\be&\al&0\\0&0&0&1\end{pmatrix}$}&3-nilpotent\\
 $a,\be\in\R$, $\al>0$. &&\\
\hline
$\begin{array}{l}
 [\bar{z},e_1]=e_2,\;[\bar{z},e_2]=-e_1.
\end{array}$&{\small$
 \begin{pmatrix}
  0&\al^{-1}&0&0\\\al^{-1}&0&\be&\ga\\0&\be&\al^2&0\\0&\ga&0&\al^2\end{pmatrix}$}&2-solvable\\
  $\be,\ga\in\R$,
$\al\not=0$.&&\\
\hline
$\begin{array}{l}
 [\bar{z},e_1]=z.
\end{array}$&{\small$
 \begin{pmatrix}
  0&\al&0&0&0\\\al&0&0&0&0\\0&0&1&a&b\\0&0&a&1+a^2&ab\\0&0&b&ab&1+b^2\end{pmatrix}$}&2-nilpotent\\
  $a,b\in\R$,
$\al>0$.&&\\
\hline
$\begin{array}{l}
 \;[\bar{z},e_1]=az,\;[\bar{z},e_2]=bz,\\\;[\bar{z},e_3]=e_1,
\;[e_1,e_3]=-z,\\\;[e_2,e_3]=-cz.
\end{array}$&{\small$
 \begin{pmatrix}
0&\al&0&0&0\\\al&0&\be&0&0\\0&\be&\al&c\al&0\\0&0&c\al&1&0\\0&0&0&0&1
\end{pmatrix}$}&3-nilpotent\\
$\al>0$, $a,b,c,\be,\in\R$.&$\al c^2\not=1.$&\\
\hline
$\begin{array}{l}
 \;[\bar{z},e_1]=e_2,\;[\bar{z},e_2]=-e_1,\\\;[\bar{z},e_3]=az.
\end{array}$&{\small$
 \begin{pmatrix}
0&\al^{-1}&0&0&0\\\al^{-1}&0&\be&\ga&0\\0&\be&\al^2&0&0
\\0&\ga&0&\al^{2}&0\\0&0&0&0&1\end{pmatrix}$}&2-solvable\\
$\al>0$, $a,\be,\ga\in\R$.&&\\
\hline
$\begin{array}{l}
 \;[\bar{z},e_1]=ae_2+be_3+cz,\\\;[\bar{z},e_2]=de_3,\;[\bar{z},e_3]=-de_2,\\
 \;[e_1,e_2]=e_3,\;[e_1,e_3]=-e_2.
\end{array}$&{\small$
 \begin{pmatrix}
  0&1&0&0&0\\1&0&0&-\al b&\al a\\0&0&\al^{-1}&0&0\\0&-\al b&0&\al &0
  \\0&\al a&0&0&\al \end{pmatrix}$}&2-solvable\\
  $(a,b,c,d)\not=0$, $\al>0$.&&\\
\hline
$\begin{array}{l}
 \;[\bar{z},e_1]=z.
\end{array}$&{\small$
 \begin{pmatrix}
     0&\al&0&0&0&0\\\al&0&0&0&0&0\\0&0&1&a&b&c\\0&0&a&1+a^2&ab&ac\\0&0&b&ab&1+b^2&bc
     \\0&0&c&ac&bc&1+c^2
    \end{pmatrix}$}&2-nilpotent\\
    $\al>0$, $a,b,c\in\R$.&&\\
\hline
\end{tabular} 

\begin{tabular}{|l|l|l|}
 \hline
$\begin{array}{l}
 \;[\bar{z},e_1]=e_2,\;[\bar{z},e_2]=-e_1,\\\;[\bar{z},e_3]=ae_4,\;[\bar{z},e_4]=-ae_3.
\end{array}$&{\small$
 \begin{pmatrix}
   0&\al^{-1}&0&0&0&0\\\al^{-1}&0&\be&\ga&\mu&\nu\\0&\be&\al^2&0&0&0
   \\0&\ga&0&\al^2&0&0
\\0&\mu&0&0&\al^2 a^2&0\\0&\nu&0&0&0&\al^2 a^2
    \end{pmatrix}$}&2-solvable \\
   $\al\not=0$, $a\not=0$, $\be,\ga,\mu,\nu\in\R$.&&\\
\hline
$\begin{array}{l}
 \;[\bar{z},e_1]=e_2,\;[\bar{z},e_2]=-e_1,\\\;[\bar{z},e_3]=az,\;[\bar{z},e_4]=bz.
\end{array}$&{\small$
 \begin{pmatrix}
   0&\al^{-1}&0&0&0&0\\\al^{-1}&0&\be&\ga&0&0\\0&\be&\al^2&0&0&0
   \\0&\ga&0&\al^2&0&0
\\0&0&0&0&1&0\\0&0&0&0&0&1
    \end{pmatrix}$}& 2-solvable\\
   $\al\not=0$, $a,b,\be,\ga\in\R$.&&\\
\hline
$\begin{array}{l}
 \;[\bar{z},e_1]=az,\;;[\bar{z},e_2]=bz,\\\;[\bar{z},e_3]=e_1,\;[\bar{z},e_4]=e_2,\\
 \;[e_1,e_3]=-\rho_1^2z,\;[e_2,e_4]=-\rho_2^2z,\\
 \;[e_1,e_4]=-\rho_1\rho_2\cos(\om_1-\om_2)z,\\
\; [e_2,e_3]=-\rho_1\rho_2\cos(\om_1-\om_2)z.
\end{array}$&{\small$
 \begin{pmatrix}
     0&1&0&0&0&0\\1&0&\be&\ga&0&0\\0&\be&\rho_1^2&\mu&0&0
     \\0&\ga&\mu&\rho_2^2&0&0\\0&0&0&0&1&0\\0&0&0&0&0&1
    \end{pmatrix}$}& 2-solvable \\
   $a,b,\be,\ga\in\R$, $\rho_1,\rho_2>0$. &$\mu=\rho_1\rho_2\cos(\om_1-\om_2)$.&\\
\hline
$\begin{array}{l}
 \;[\bar{z},e_1]=az,\;;[\bar{z},e_2]=bz,\\\;[\bar{z},e_3]=cz,\;[d,e_4]=e_1,\\
 \;[e_1,e_4]=-z,\;[e_2,e_4]=-d z,\\\;[e_3,e_4]=-e z.
\end{array}$&{\small$
 \begin{pmatrix}
     0&\al&0&0&0&0\\\al&0&\be&0&0&0\\0&0&\al&\al d&\al e&0\\0&0&\al
d&1&0&0\\0&0&\al e&0&1&0\\0&0&0&0&0&1
    \end{pmatrix}$}& 3-nilpotent\\
    $\al>0$, $a,b,c,d,e,\be\in\R$.&$\al(d^2+e^2)\not=1$.&\\
\hline
$\begin{array}{l}
 \;[\bar{z},e_1]=e_2,\;;[\bar{z},e_2]=-e_1,\\\;[\bar{z},e_3]=az,\;[d,e_4]=e_3,\\
 \;[e_3,e_4]=-z.
\end{array}$&{\small$
 \begin{pmatrix}
     0&\al&0&0&0&0\\\al&0&\ga&\mu&\nu&0\\0&\ga&\be&0&0&0\\
     0&\mu&0&\be&0&0\\0&\nu&0&0&\al&0\\0&0&0&0&0&\be
    \end{pmatrix}$}& 2-solvable\\
    $\al,\be>0$, $a,\ga,\mu,\nu\in\R$.&&\\
\hline
{$\begin{array}{l}
\;[{}{\bar{z}},{}{e}_1]={}{e}_2,\; \;[{}{\bar{z}},{}{e}_2]=-{}{e}_1,\\
  \;[{}{\bar{z}},{}{e}_3]=Y{}{e}_2+X{}{e}_1+{}{e}_4,\\
  \;[{}{\bar{z}},{}{e}_4]=X{}{e}_2-Y{}{e}_1-{}{e}_3,\\
\;[{}{e}_1,{}{e}_3]=-X{}{z},\;[{}{e}_1,{}{e}_4]=
Y{}{z},\\\;[{}{e}_2,{}{e}_3]=-Y{}{z},
\;[{}{e}_2,{}{e}_4]=-X{}{z},\\\;[{}{e}_3,{}{e}_4]=2\rho^2{z}.
\end{array}$}&{\footnotesize$ \begin{pmatrix}
   0&1&0&0&0&0\\
   1&0&\al&\be&\ga&\mu\\
   0&\al&1&0&Y&X\\
   0&\be&0&1&-X&Y\\
   0&\ga&Y&-X&1+\rho^2&0\\
   0&\mu&X&Y&0&1+\rho^2
  \end{pmatrix}$}& 3-solvable\\
 $\rho>0$, $\al,\be,\ga,\mu,\om\in\R$. &$X=\rho\cos\om,\;Y=\rho\sin\om.$&\\
\hline
$\begin{array}{l}
 \;[\bar{z},e_1]=ae_2+be_3+ce_4+dz,\\\;[\bar{z},e_2]=ez,\;[\bar{z},e_3]=fe_4,\\\;
 [\bar{z},e_4]=-fe_3,
 \;[e_1,e_2]=az,\\\;[e_1,e_3]=e_4,\;[e_1,e_4]=-e_3
\end{array}$&{\small$
 \begin{pmatrix}
   0&\al&0&0&0&0\\
   \al&0&0&0&-c&b\\
   0&0&\be&\ga&0&0\\
   0&0&\ga&\al&0&0\\
   0&-c&0&0&1&0\\
   0&b&0&0&0&1
  \end{pmatrix}$}& 2-solvable\\
  $\al,\be>0$, $a,b,c,d,e,f,\ga\in\R$.&$\al\be-\ga^2=1.$&\\
\hline
\end{tabular}

Malika AitBenhaddou\\               
Facult\'e des Sciences Meknes\\ 
Meknes Morocco\\
maitbenhaddou@gmail.com\\\\
Mohamed Boucetta\\
Universit\'e Cadi-Ayyad\\
Facult\'e des Sciences Gueliz\\
Marrakech (Morocco)\\
mboucetta2@yahoo.fr\\\\
Hicham Lebzioui\\
Facult\'e des Sciences Meknes\\
Meknes Morocco\\
hlebzioui@gmail.com

\end{document}